\documentclass[a4,11pt]{article}
\usepackage{amsmath,amsfonts,amsthm,latexsym,amssymb,mathrsfs}
\usepackage{bbm}

\begin{document}
\pagestyle{myheadings} \markboth{  }{ \hskip6truecm \rm Enrico  Boasso}

\title{ Characterizations of Fredholm Pairs \\and Chains in Hilbert Spaces}

\author{ Enrico Boasso} 
\date{ }
\maketitle
\begin{abstract}{In this work characterizations of Fredholm pairs and chains of Hilbert space operators are given. Following a well-known idea of several variable operator
theory in Hilbert spaces, the aforementioned objects are characterized in terms of Fredholm linear and bounded maps. Furthermore, as an application of the main results of this work, direct proofs of the stability properties of Fredholm pairs and chains in Hilbert spaces are obtained.\par
\vskip.2cm
\noindent \it AMS 2000 Subject Classification\rm: Primary 47A13; Secondary
47A53, 47A55.\par
\vskip.2cm
\noindent \it Key words and phrases\rm: Fredholm pairs and chains, Fredholm 
self-adjoint operators, index, stability properties.}
\end{abstract}


\baselineskip=16truept
\vskip.5cm
\centerline{\bf 1. INTRODUCTION}\par
\vskip.5cm

\indent In multiparameter operator theory there are objects that in the frame of
Hilbert spaces can be described in terms of suitable linear and continuous 
maps. It is important to remark that this reduction of the complexity of certain
problems from several variable operator theory to the classical setting of one single operator allows not only to recover properties and techniques from one dimensional operator theory, but also to simplify proofs. For instance, as regard Hilbert space complexes, the exactness (resp. Fredholmness)
of such an object can be determined by the invertibility (resp. the Fredholmness)
of a Hilbert space operator, see [11], [5] and [12, Chap. III]. Furthermore, this connection
between one and several variable operator theory was for
the first time developed to give a characterization of the Taylor joint spectrum in Hilbert space in terms of the invertibility of a single Hilbert space linear and bounded map, see [8] and [9].
In addition, the application of this idea led to many new results
in the area under consideration, see for example [8], [9], [3], [4], [10], [11], [5] and [12].  
\par

\indent On the other hand, Fredholm pairs were studied in the works [1] and [2],  
where the main stability properties of such objects were also proved. Roughly speaking, the aforementioned pairs consist in an extension of the notion of Fredholm operator to multiparameter
spectral theory, which is closely related to the concept of Fredholm Banach space complex, see [1] and [2].
However, Fredholm pairs have not been studied in the frame of Hilbert
spaces yet. In this work two characterization of these objects are given.
In fact, following the idea mentioned in the first paragraph, Fredholm pairs are characterized in terms of Fredholm 
Hilbert space linear and bounded maps. Moreover, these characterizations will be
applied to directly prove the stability properties of the objects under consideration, as well as to study dual Fredholm pairs, see section 4.\par

\indent A generalization of the concept of Fredholm Banach 
spaces complex is the notion of Fredholm chain, which is closely related
to the one of Fredholm pair, see [6]. In this work, thanks to the above-mentioned results on Fredholm pairs, two characterizations of Fredholm chains in Hilbert spaces will be given. In fact, these objects will be characterized in terms of Fredholm Hilbert space operators. Furthermore, the stability properties of the objects under consideration will be proved and dual Fredholm chains
will be studied, see section 6.
\par

\indent The article is organized as follows. In the next section some definitions and facts needed for the present work are reviewed. In sections 3 and 4 the main results of this work are proved. In addition, in section 4 dual Fredholm pairs are also studied. In section 5, as an application of the characterization studied in section 4, the stability properties
of Fredholm pairs in Hilbert spaces are proved. Finally, in section 6, Fredholm chains of Hilbert space operators
are considered. In fact, these objects are characterized, their stability properties are
proved, and dual Fredholm chains are studied.\par
\vskip.5cm
\centerline{\bf 2. PRELIMINARY DEFINITIONS AND FACTS}\par
\vskip.5cm
\indent Since all the operators considered in this article will be defined on Hilbert spaces, all the definitions and facts
reviewed will be restricted to this class of spaces and maps. For a general presentation, see the works [1], [2] and [6].\par
 
 \indent From now on, $H_1$ and $H_2$ denote two Hilbert spaces, $L(H_1,H_2)$ the algebra 
of all linear and continuous maps defined on $H_1$ with values in $H_2$, and $K(H_1,H_2)$ the closed
ideal of all compact operators of $L(H_1,H_2)$. As usual, when $H_1=H=H_2$, $L(H,H)$ and $K(H,H)$ are denoted
by $L(H)$ and $K(H)$ respectively. For
every $S\in L(H_1,H_2)$, the null space of $S$ is denoted by
$N(S)=\{x\in H_1\colon \hbox{  }S(x)=0\}$, and the range of $S$ by $R(S)=\{ y\in H_2: \hbox{  }\exists
\hbox{ } x\in H_1 \hbox{ such that }y=S(x)\}$. Next follows the definition of Fredholm pair, 
see for instance [1].\par

\indent {\it Definition\rm} 2.1. Let $H_1$ and $H_2$ be two Hilbert spaces. Let $S\in L(H_1,H_2)$ and
$T\in L(H_2,H_1)$ be such that the following dimensions are finite:
$$
a\colon =\dim N(S)/(N(S)\cap R(T)), \hbox{ }b\colon = \dim R(T)/(N(S)\cap R(T)),
$$
$$
c\colon = \dim N(T)/(N(T)\cap R(S)), \hbox{ }d\colon =\dim R(S)/(N(T)\cap R(S)).
$$
A pair $(S,T)$ with the above properties is called a Fredholm pair (see [1]).\par
\indent Let $P(H_1,H_2)$ denote the set of all Fredholm pairs. If $(S,T)\in P(H_1,H_2)$, then
the index of $(S,T)$ is defined by the equality
$$
\hbox{\rm ind} \hbox{ } (S,T)\colon = a -b -c +d.
$$

\indent Before going on, several properties of Fredholm pairs are recalled, see [1].\par

\indent {\it Remark \rm } 2.2. First of all, observe that if $S\in L(H_1,H_2)$ is a Fredholm operator, then 
$(S,0)$ is a Fredholm pair. Furthermore, $\hbox{\rm ind}\hbox{ } S=\hbox{\rm ind}\hbox{ } (S,0)$. Consequently, the definition of Fredholm pair extends the notion of Fredholm operator to several variable operator theory.\par
\indent In second place, note that if $(S,T)\in P(H_1,H_2)$, then $(T,S)\in P(H_2,H_1)$ and
$$
\hbox{\rm ind}\hbox{ } (T,S)=-\hbox{\rm ind}\hbox{ } (S,T).
$$ 
\indent Finally, if $(S,T)\in P(H_1,H_2)$, then $R(S)$ and $N(T) +R(S)$ are closed 
subspaces in $H_2$. Similarly, $R(T)$ and $N(S)+R(T)$ are closed subspaces in $H_1$.\par

\indent Next follows the definition of Fredholm chains, see for instance [6].

\indent {\it Definition \rm} 2.3. A Fredholm chain $(H,\delta)$ is 
a sequence of spaces and maps
$$
0\to H_n\xrightarrow{\delta_n} H_{n-1}\to\ldots\to H_1\xrightarrow{\delta_1} H_0\to 0,
$$
where $H_p$ are Hilbert spaces, and $\delta_p\in L(H_p,H_{p-1})$ are bounded operators such that
$$
N(\delta_p)/(N(\delta_p)\cap R(\delta_{p+1}))\hbox{ and }
R(\delta_{p+1})/(N(\delta_p)\cap R(\delta_{p+1}))
$$   
are finite dimensional subspaces of $H_p$, p = 0, $\ldots$, n.
Formally, it is assumed that $H_p=0$ and
$\delta_p=0$, for p $<$ 0 and p $\ge$ n+1.\par

\indent Given a Fredholm chain, it is possible
to associate to it an index. In fact, if $(H,\delta)$ is such an object, then define 

\begin{align*}
\hbox{ind } (H,\delta)=\sum_{p=0}^n(-&1)^p (\dim\hbox{ }
N(\delta_p)/(N(\delta_p)\cap R(\delta_{p+1}))\\
&-\dim R(\delta_{p+1})/(N(\delta_p)\cap R(\delta_{p+1}))),\\
\end{align*}
see [6].\par 
\indent Recall that in [6] it was introduced the more general concept of 
semi-Fredholm chain.
However, since the main concern of this article consists in Fredholm objects,
only Fredholm chains will be considered. Furthermore,
observe that since
$\dim R(\delta_{p-1}\delta_p)$ are finite dimensional, p = 1, $\ldots$, n, a Fredholm
chain $(H,\delta)$
is a particular case of what in [7] was called an  essential complex of Banach spaces.\par
\indent In the following remark, the relationship between Fredholm pairs
and chains is considered.\par
\indent {\it Remark \rm } 2.4. Consider $(H,\delta)$ a sequence of spaces and maps 
$$
0\to H_n\xrightarrow{\delta_n} H_{n-1}\to\ldots\to H_1\xrightarrow{\delta_1} H_0\to 0,
$$
where $H_p$ are Hilbert spaces, and $\delta_p\in L(H_p,H_{p-1})$ are bounded operators. In addition, assume that $H_p=0$ and $\delta_p=0$, for
p $<$ 0 and p $\ge$ n+1.\par
\indent Next, associate to this sequences the Hilbert spaces 
$$
H_1=\bigoplus_{p=2k}H_p,\hskip1cm H_2=\bigoplus_{p=2k+1}H_p,
$$
and the maps $S\in L(H_1,H_2)$ and $T\in L(H_2,H_1)$ defined as 
$$
S=\bigoplus_{p=2k}\hbox{ } \delta_p,\hskip1cm T=\bigoplus_{p=2k+1}\hbox{ } \delta_p,
$$
where $H_p=0$ and $\delta_k=0$, when p $<$ 0 and p $\ge$ n+1.\par

\indent Since
$$
R(ST)=\bigoplus_{p=2k} R(\delta_{p}\delta_{p+1})\in L(H_2),\hskip.5cm
R(TS)=\bigoplus_{p=2k+1}R(\delta_{p}\delta_{p+1})\in L(H_1),
$$
dim $R(ST)$ and dim $R(TS)$ are finite dimensional if and only if
dim $R(\delta_p\delta_{p+1}) =$ dim $R(\delta_{p+1})/(N(\delta_p)\cap R(\delta_{p+1}))$ are finite dimensional, p = 0, $\ldots$, n.\par
\indent Furthermore, since it is clear that
$$
N(S)/(N(S)\cap R(T))=\bigoplus_{p=2k}N(\delta_p)/(N(\delta_p)\cap R(\delta_{p+1})), 
$$
$$
N(T)/(N(T)\cap R(S))=\bigoplus_{p=2k+1}N(\delta_p)/(N(\delta_p)\cap R(\delta_{p+1})),
$$
the sequence $(H,\delta)$ is a Fredholm chain if and only if $(S,T)$ is a Fredholm
pair.\par
\indent Finally, in this case a straightforward calculation shows that
$$
\hbox{\rm ind } (H,\delta)=\hbox{\rm ind } (S,T).
$$

\centerline{\bf 3. THE FIRST CHARACTERIZATION}
\vskip.5cm
\indent In this section, the first characterization of Fredholm pairs of Hilbert space
operators is given. 
To this end, the argument develops 
ideas of [12, Chap. III, Crollary 7.4], see also [8, Theorem 2.1] and [5, Proposition 2.1].\par 
\indent In order to prove the main result of this section, several previous propositions are
needed. In first place, some preparation is presented.\par

\indent Let $H_1$ and $H_2$ be two Hilbert spaces and consider 
$S\in L(H_1,H_2)$ and $T\in L(H_2,H_1)$, two operators such that
$\dim\hbox{ } R(ST)$ and $\dim\hbox{ } R(TS)$ are finite dimensional. Define
the Hilbert space $H$ as the orthogonal direct sum of $H_i$, i=1, 2,
that is $H=H_1\oplus H_2$. Now well, if $U\in L(H)$ is the linear
and continuous map

$$U=\begin{pmatrix}
0&T\\
S&0\\\end{pmatrix} ,
$$

then consider the self-adjoint operator $V=U+U^*\in L(H)$. The following proposition is the first step to 
the first characterization of Fredholm pairs in Hilbert spaces.\par
 
\indent PROPOSITION 3.1. \it Let $H_1$ and $H_2$ be two Hilbert spaces and consider
two bounded and linear maps $S\in L(H_1,H_2)$ and $T\in L(H_2,H_1)$ such that 
$R(TS)$ and $R(ST)$ are finite dimensional subspaces of $H_1$ and
$H_2$ respectively. Then, if $U$ and $V$ are
the operators defined above, the following statements are equivalent:\par
\hskip.2cm i) $V=U+U^*$ is a Fredholm operator,\par
\hskip.2cm ii) $V^2$ is a Fredholm operator,\par
\hskip.2cm iii) $UU^*+U^*U$ is a Fredholm operator,\par
\hskip.2cm iv) $TT^*+S^*S\in L(H_1)$ and $SS^*+T^*T\in L(H_2)$ are
Fredholm operators.\rm\par
\indent {\it Proof. \rm } First of all, it is clear that $V$ is a Fredholm operator if and only if $V^2$ is. 
In addition, note that $V^2=U^2 + U^{*2}+UU^*+U^*U$, where
$$
U^2=\begin{pmatrix}
TS&0\\
0&ST\\
\end{pmatrix},\hbox{  }
U^{*2}=\begin{pmatrix}
(ST)^*&0\\
0&(TS)^*\\
\end{pmatrix},
$$
and
$$UU^*+U^*U=\begin{pmatrix}
TT^*+S^*S&0\\
0&SS^*+T^*T\\
\end{pmatrix}.
$$ 
\indent Now well, since dim $R(ST)$ and dim $R(TS)$ are finite dimensional,
dim $R(U^2)$ and dim $R(U^{*2})$ also are finite dimensional. Therefore,
$V^2$ is a Fredholm operator if and only if 
$UU^* +U^*U$ is a Fredholm operator, which is equivalent to the fact that
$TT^*+S^*S$ and $SS^*+T^*T$ are Fredholm operators.\qed\par
\indent The following remarks are needed to prove the main result of the 
present section.\par
\indent {\it Remark \rm } 3.2. Let $H$ be a Hilbert space and
consider $T\in L(H)$ such that $T^2=0$, that is $R(T)\subseteq N(T)$. Then, according to [12, Chap. III, Lemma 7.3],
$R(T)=N(T)$ if and only if the self-adjoint operator $T+T^*$ is 
invertible on $H$.
Now well, suppose that $R(T)$ is closed. Then, according to the argument 
in the aforementioned Lemma, it is easy to prove that dim $N(T)/R(T)$ is finite dimensional 
if and only if $T+T^*$ is a Fredholm self-adjoint operator.\par

\indent {\it Remark \rm } 3.3. Let $H_1$ and $H_2$ be two Hilbert
spaces and consider two bounded and linear maps $S\in L(H_1,H_2)$ and $T\in L(H_2,H_1)$ such that     
$F_1=R(TS)$ and $F_2=R(ST)$ are finite dimensional subspaces of $H_1$ and $H_2$ respectively. Decompose $H_1$ and $H_2$ as the following orthogonal direct sum:
$$
H_1=F_1\oplus \mathcal{H}_1,\hskip1cm H_2=F_2\oplus \mathcal{H}_2,
$$
where $\mathcal{H}_i=F_i^{\perp}$, i=1, 2.\par
\indent  Next consider the operators
$\mathcal{S}\in L(\mathcal{H}_1,\mathcal{H}_2)$ and $\mathcal{T}\in L(\mathcal{H}_2,\mathcal{H}_1)$ defined by
$$
\mathcal{S}= P_2\circ (S\mid \mathcal{H}_1),\hskip1cm \mathcal{T}= P_1\circ (T\mid \mathcal{H}_2),
$$
where $P_i\colon H_i\to\mathcal{H}_i$ is the orthogonal projection onto $\mathcal{H}_i$, i=1, 2. It is clear that 
$$
\mathcal{S}\circ\mathcal{T}=0,\hskip1cm \mathcal{T}\circ\mathcal{S}=0.
$$
\indent Now well, according to [1, Remark 2.1], it is easy to prove that $(S,T)$ is a Fredholm pair if and only if $(\mathcal{S},\mathcal{T})$ is a Fredholm pair, which is equivalent to the fact that dim $N(\mathcal{S})/ R(\mathcal{T})$ and dim $N(\mathcal{T})/R(\mathcal{S})$ are finite
dimensional.\par 
\indent Furthermore, in this case
$$
\hbox{\rm ind } (S,T)=\hbox{\rm ind} (\mathcal{S},\mathcal{T}) - \hbox{\rm dim } R(ST) +\hbox{ \rm dim } R(TS).
$$  
\indent The next result is the key step for the first characterization of 
Fredholm pairs in Hilbert spaces.\par

\indent PROPOSITION 3.4.  \it Let $H_1$ and $H_2$ be two Hilbert spaces and consider
two bounded and linear maps $S\in L(H_1,H_2)$ and $T\in L(H_2,H_1)$ such that     
$R(TS)$ and $R(ST)$ are finite dimensional subspaces of $H_1$ and $H_2$
respectively. Then, $(S,T)$ belongs to $P(H_1,H_2)$ if and only if  $TT^*+S^*S\in L(H_1)$ and $SS^*+T^*T\in L(H_2)$ are Fredholm operators.\par \rm
\indent {\it Proof. \rm}First of all, consider the finite dimensional subspaces $F_1= R(TS)$ and 
$F_2=R(ST)$, and decompose $H_1$ and $H_2$ as in Remark 3.3, that is
$$
H_1=F_1\oplus \mathcal{H}_1,\hskip1cm H_2=F_2\oplus \mathcal{H}_2,
$$
where $\mathcal{H}_i=F_i^{\perp}$, i=1, 2.\par
\indent In addition, consider the operators $\mathcal{S}\in L(\mathcal{H}_1,\mathcal{H}_2)$ and 
$\mathcal{T}\in L(\mathcal{H}_2,\mathcal{H}_1)$ defined in Remark 3.3.\par
\indent Now suppose that $(S,T)$ is a Fredholm pair. Then,
according to Remark 3.3 or to [1, Remark 2.1], $(\mathcal{S},\mathcal{T})$ is a Fredholm pair, that is
dim $N(\mathcal{S})/ R(\mathcal{T})$ and dim $N(\mathcal{T})/R(\mathcal{S})$ are finite
dimensional.\par 
\indent Now well, if $\mathcal{H}=\mathcal{H}_1\oplus \mathcal{H}_2$,
and if 
$$
\mathcal{U}=\begin{pmatrix}
0&\mathcal{T}\\
\mathcal{S}& 0\\
\end{pmatrix},
$$
then it is easy to prove that $R(\mathcal{U})$ is a closed subspace of $\mathcal{H}$, ${\mathcal{U}}^2=0$, and  dim $N(\mathcal{U})/R(\mathcal{U})$ is finite dimensional. Consequently, according to
Remark 3.2, $\mathcal{V}=\mathcal{U}+{\mathcal{U}}^*$ is also a Fredholm operator,
and according to Proposition 3.1, $\mathcal{T}{\mathcal{T}}^*+{\mathcal{S}}^*\mathcal{S}\in L(\mathcal{H}_1)$
and $\mathcal{S}{\mathcal{S}}^*+{\mathcal{T}}^*\mathcal{T}\in L(\mathcal{H}_2)$ are Fredholm operators.\par 

\indent On the other hand, if $\mathcal{S}$ (resp. $\mathcal{T}$) is extended to $H_1$ (resp. $H_2$)
by setting $\mathcal{S}\mid F_1\equiv 0$ (resp. $\mathcal{T}\mid F_2\equiv 0$),
then it is clear that there are operators $S_1\in L(H_1,H_2)$
and $T_1\in L(H_2,H_1)$ such that $R(S_1)$ and $R(T_1)$ are finite dimensional
and 
$$
S=\mathcal{S}+S_1,\hskip1cm T=\mathcal{T}+T_1,
$$   
where $\mathcal{S}\in L(H_1,H_2)$ and $\mathcal{T}\in L(H_2,H_1)$ also denote
the extension of $\mathcal{S}\in L(\mathcal{H}_1,\mathcal{H}_2)$ and $\mathcal{T}\in L(\mathcal{H}_2,\mathcal{H}_1)$ respectively.\par
\indent Now well, a straightforward calculation proves that there are two operators
$K_1\in L(H_1)$ and $K_2\in L(H_2)$ whose ranges are finite dimensional and such 
that
$$
TT^*+S^*S=\mathcal{T}{\mathcal{T}}^*+{\mathcal{S}}^*\mathcal{S}
+K_1,\hbox{  } SS^*+T^*T=\mathcal{S}{\mathcal{S}}^*+{\mathcal{T}}^*\mathcal{T}
+K_2.
$$
Therefore, $TT^*+S^*S$ and $SS^*+T^*T$ are Fredholm operators.\par 

\indent Conversely, according to the previous argument, if $TT^*+S^*S$ and $SS^*+T^*T$ are Fredholm operators, then it is clear that
$\mathcal{T}{\mathcal{T}}^*+{\mathcal{S}}^*\mathcal{S}\in L(\mathcal{H}_1)$ and $\mathcal{S}{\mathcal{S}}^*+{\mathcal{T}}^*\mathcal{T}\in L(\mathcal{H}_2)$
are Fredholm opeartors. Furthermore, according to Proposition 3.1, $\mathcal{V}=\mathcal{U}+
\mathcal{U}^*\in L(\mathcal{H})$ is a Fredholm operator, and according to Remark 3.2,
dim $N(\mathcal{U})/R(\mathcal{U})$ is finite dimensional.     \par 
\indent Now well, 
$$
N(\mathcal{U})=N(\mathcal{S})\oplus N(\mathcal{T}),\hskip.5cm R(\mathcal{U})=R(\mathcal{T})\oplus R(\mathcal{S}).
$$
Consequently, $(\mathcal{S},\mathcal{T})$ is a Fredholm pair. However, since $R(ST)$ and $R(TS)$ are finite dimensional 
subspaces of $H_2$ and $H_1$ respectively, according to Remark 3.3 or to [1, Remark 2.1], $(S,T)$ is a Fredholm pair.\qed\par

\indent In the next theorem the first characterization of Fredholm pairs
in Hilbert spaces is presented.\par
\indent THEOREM 3.5. \it Let $H_1$ and $H_2$ be two Hilbert spaces and consider
two bounded and linear maps $S\in L(H_1,H_2)$ and $T\in L(H_2,H_1)$ such that     
$R(TS)$ and $R(ST)$ are finite dimensional subspaces of $H_1$ and $H_2$
respectively. Then, with the notations of Proposition 3.1, the following statements are equivalent:\par
\hskip.2cm i) $(S,T)$ is a Fredholm pair\par
\hskip.2cm ii) $V=U+U^*$ is a Fredholm operator,\par
\hskip.2cm iii) $V^2$ is a Fredholm operator,\par
\hskip.2cm iv) $UU^* + U^*U$ is a Fredholm operator,\par
\hskip.2cm v) $TT^*+S^*S\in L(H_1)$ and $SS^*+T^*T\in L(H_2)$ are
Fredholm operators.\par \rm
\indent {\it Proof. \rm } It is a consequence of Propositions 3.1 and 3.4.\qed\par  
\vskip.5cm
\centerline{\bf 4. THE SECOND CHARACTERIZATION}
\vskip.5cm

\indent In Theorem 3.5 the condition of being a Fredholm pair is expressed in terms of Fredholm self-adjoint operators. This formulation has the disadvantage that the index of a Fredholm pair can not be related to the index of any of
the linear and contimuos maps considered in the aforementioned theorem. In this section, however, it is proved another characterization of Fredholm pairs in Hilbert spaces which precisely has the
advantage that the index of a Fredholm pair is expressed in terms of
the index of a Fredholm operator defined between two Hilbert spaces.
In addition, as a first appilcation of this characterization, dual Fredholm 
pairs will be studied. On the other hand, the argument of the main result of this section develops ideas of
[11, Theorem 1.2] and [12, Chap. III, Theorem 7.1].\par
\indent First of all, in order to
prove the second characterization of Fredholm pairs in Hilbert spaces, some preparation is needed.\par

\indent {PROPOSITION 4.1.} \it  Let $H_1$ and $H_2$ be
two Hilbert spaces and consider two operators $S\in L(H_1,H_2)$
and $T\in L(H_2,H_1)$ such that $T S=0$ and $S T=0$. Then, the following
statements are equivalent:\par
\hskip.2cm i) dim $N(S)/R(T)$ and dim $N(T)/R(S)$ are finite dimensional,\par
\hskip.2cm ii) $S+T^*\in L(H_1,H_2)$ is a Fredholm operator,\par
\hskip.2cm iii)  $T+S^*\in L(H_2,H_1)$ is a Fredholm operator.\par
\indent Furthermore, in this case $(S,T)$ is a Fredholm pair and
$$
\hbox{ \rm ind}\hbox{ } (S,T)=\hbox{ }\hbox{\rm ind}\hbox{  } (S+T^*)=- \hbox{ }\hbox{\rm ind}\hbox{ }  (T+S^*).
$$\rm
\indent{\it Proof. \rm } First of all, it is clear that $S+T^*$ is a Fredholm operator
if and only if $T+S^*$ is, and in this case ind $(T+S^*)$ = - ind $(S+T^*)$.\par
\indent Next, note that if dim $N(S)/R(T)$ and  dim $N(T)/R(S)$
are finite dimensional, then $R(S)$ and $R(T)$ are closed subspaces
of $H_2$ and $H_1$ respectively. Conversely, if $S+T^*$ is a 
Fredholm operator, then $R(S)$ and $R(T)$ are closed.\par
\indent In fact, consider $(x_n)_{n\in\Bbb N}\subseteq N(S)^{\perp}$
such that $(S(x_n))_{n\in\Bbb N}$ converges to $v\in H_2$. Since
the closure of $R(T)$ is contained in $N(S)$, 
$$
N(S)^{\perp}\subseteq \overline{R(T)}^{\perp}=R(T)^{\perp}\subseteq N(T^*).
$$   
\indent In particular, $(S+T^*)(x_n)= S(x_n)$, and $v\in R(S+T^*)$. Thus,
there is $x\in H_1$ such that $(S+T^*)(x)=v$. However, $(S(x_n-x))_{n\in\Bbb N}$
converges to $T^*(x)$. Then,
$$
T^*(x)\in R(T^*)\cap \overline{R(S)}\subseteq R(T^*)\cap N(T)\subseteq R(T^*)\cap R(T^*)^{\perp}=0.
$$
\indent Consequently, $v=S(x)\in R(S)$.\par
\indent Since $T+S^*$ is Fredholm if and only if $T+S^*$ is, a similar argument proves that $R(T)$ is a closed
subspace of $H_1$.\par
\indent Now well, since it has been proved that in the conditions of the proposition
$R(S)$ and $R(T)$ are always closed subspaces of $H_2$ and $H_1$ respectively, in the rest
of the proof it will be assumed this property for $S$ and $T$.\par

\indent Decompose, then, $H_1$ and $H_2$ as the orthogonal direct
sum of the following subspaces:
$$
H_1=(R(T)\oplus N_1)\oplus L_1,\hskip.5cm
H_2=(R(S)\oplus N_2)\oplus L_2,
$$
where 
$$
R(T)\oplus N_1=N(S),\hskip.1cm L_1=N(S)^{\perp},\hskip.1cm 
R(S)\oplus N_2=N(T),\hskip.1cm L_2=N(T)^{\perp}.
$$
\indent Now well, according to the above orthogonal direct sum, and
since $N(T^*)=R(T)^{\perp}=N_1\oplus L_1$ and $R(T^*)=N(T)^{\perp}=L_2$, it is
clear that the
operator $S+T^*$ can be presented in the following matricial form:
$$ 
\begin{pmatrix} 0&0&\mathcal{S}\\
               0&0&0\\
               \mathcal{T}^*&0&0\\\end{pmatrix},
$$  
where 
$$
\mathcal{S}=S\mid_{L_1}\colon L_1\xrightarrow{\cong} R(S),\hskip.2cm
\mathcal{T}=T\mid_{L_2}\colon L_2\xrightarrow{\cong} R(T).
$$
\indent Moreover, this matricial decomposition gives  
$$
N(S+T^*)=N_1,\hskip.5cm R(S+T^*)^{\perp}=N_2.
$$
\indent Therefore, $N(S)/R(T)$ and $N(T)/R(S)$ are finite dimensional 
subspaces of $H_1$ and $H_2$ respectively if and only if $S+T^*$ is a 
Fredholm operator. \par
\indent In addition, in this case, since $ST=0$ and $TS=0$,
and since dim $N_1$ = dim $N(S)/R(T)$ and dim $N_2$ = dim $N(T)/R(S)$,
$(S,T)\in P(H_1,H_2)$ and 
$$
\hbox{\rm ind } (S,T)=\hbox{ \rm ind } (S+T^*).\qed
$$

\indent Next follows the main result of the present section.\par

\indent{THEOREM 4.2.} \it Let $H_1$ and $H_2$ be two Hilbert spaces and consider
two bounded and linear maps $S\in L(H_1,H_2)$ and $T\in L(H_2,H_1)$ such that     
$R(TS)$ and $R(ST)$ are finite dimensional subspaces of $H_1$ and $H_2$
respectively. Then, the following statements are equivalent:\par
\hskip.2cm i) $(S,T)$ is a Fredholm pair,\par
\hskip.2cm ii) $S+T^*\in L(H_1,H_2)$ is a Fredholm operator,\par
\hskip.2cm iii)  $T+S^*\in L(H_2,H_1)$ is a Fredholm operator.\par
\indent Furthermore, in this case
$$
\hbox{ } \hbox{\rm ind}\hbox{ } (S,T)=\hbox{ } \hbox{\rm ind}\hbox{  } (S+T^*)=- \hbox{  }\hbox{\rm ind}\hbox{  } (T+S^*).
$$\rm

\indent{\it Proof. \rm } First of all, consider the finite dimensional subspaces $F_1= R(TS)$ and 
$F_2=R(ST)$, and decompose $H_1$ and $H_2$ as in Remark 3.3, that is
$$
H_1=F_1\oplus \mathcal{H}_1,\hskip1cm H_2=F_2\oplus \mathcal{H}_2,
$$
where $\mathcal{H}_i=F_i^{\perp}$, i=1, 2.\par
\indent In addition, consider the operators $\mathcal{S}\in L(\mathcal{H}_1,\mathcal{H}_2)$ and 
$\mathcal{T}\in L(\mathcal{H}_2,\mathcal{H}_1)$ defined in Remark 3.3. Therefore, according to
Remark 3.3 or to [1, Remark 2.1], $(S,T)$ is a Fredholm pair if and only if 
$(\mathcal{S},\mathcal{T})\in P(\mathcal{H}_1,\mathcal{H}_2)$, equivalently, 
$N(\mathcal{S})/R(\mathcal{T})$
and $N(\mathcal{T})/ R(\mathcal{T})$ are finite dimensional subspaces of
$\mathcal{H}_1$ and $\mathcal{H}_2$ respectively. However, 
according to
Proposition 4.1, this is equivalent to the fact that $\mathcal{S} +\mathcal{T}^*\in L(\mathcal{H}_1,
\mathcal{H}_2)$ is a Freholm operator.\par

\indent Next, as in Proposition 3.4, extend
$\mathcal{S}$ (resp. $\mathcal{T}$) to $H_1$ (resp. $H_2$)
by setting $\mathcal{S}\mid F_1\equiv 0$ (resp. $\mathcal{T}\mid F_2\equiv 0$),
and denote this extension by $\mathcal{S}$ (resp. $\mathcal{T}$).
 In addition, again as in Proposition 3.4, 
consider operators $S_1\in L(H_1,H_2)$
and $T_1\in L(H_2,H_1)$ such that $R(S_1)$ and $R(T_1)$ are finite dimensional
and 
$$
S=\mathcal{S}+S_1,\hskip1cm T=\mathcal{T}+T_1.
$$ 
\indent Now well, since $F_1$ and $F_2$ are finite dimensional
subspaces of $H_1$ and $H_2$, $\mathcal{S}+\mathcal{T}^*\in L(\mathcal{H}_1,
\mathcal{H}_2)$ is a Fredholm operator if and only if 
$\mathcal{S}+\mathcal{T}^*\in L(H_1,H_2)$ is a Fredholm bounded a linear map.
However, since $R(S_1)$ and $R(T^*_1)$ are finite dimensional subspaces of
$H_2$, and since
$$
S+T^* =(\mathcal{S} +\mathcal{T}^*) + (S_1+T_1^*),
$$
$\mathcal{S} +\mathcal{T}^*\in L(H_1,H_2)$ is a Fredholm operator if and only if
$S+T^*$ is a Fredholm linear and continuous map.
Moreover, in this case
$$
\hbox{\rm ind } (S+T^*) =\hbox{ \rm ind }(\mathcal{S} +\mathcal{T}^*).
$$
\indent However, according to Proposition 4.1, for $\mathcal{S}+\mathcal{T}^*\in L(H_1,H_2)$, 
$$
\hbox{ \rm ind } (\mathcal{S}+\mathcal{T}^*)= \hbox{ \rm ind } (\mathcal{S},\mathcal{T}) -\hbox{ \rm dim } R(ST)
+\hbox{ \rm dim } R(TS),
$$
and according to Remark 3.3 or to [1, Remark 2.1],
$$
\hbox{\rm ind } (S+T^*)=\hbox{ \rm ind } (S,T).
$$
\indent Finally, it is clear that $S+T^*$ is a Fredholm operator if and only if
$T+S^*$ is, and in this case, ind $(S+T^*)$ = $-$ ind $(T+S^*)$.\qed\par  

\indent As a first application of Theorem 4.2, the relationship between Fredholm
pairs and adjoint operators in Hilbert spaces is studied. To this end, some
preparation is needed.\par

\indent Let $H_i$, i=1, 2, be two Hilbert spaces, $S$ belong to $L(H_1,H_2)$ and $T$ to $L(H_2,H_1)$. Denote $(S,T)^*\equiv (T^*,S^*)$, where $S^*$ (resp. $T^*$) is the adjoint map of $S$ (resp. $T$). Note that $(S,T)^{**}=(T^*,S^*)^*=(S,T)$.\par

\indent { THEOREM 4.3.} \it Let $H_1$ and $H_2$ be two Hilbert spaces and consider 
$(S,T)\in L(H_1,H_2)\times L(H_2,H_1)$. Then, $(S,T)$ is a Fredholm pair if and 
only if $(T^*, S^*)$ is. Furthermore, in this case
$$
\hbox{ }\hbox{\rm ind}\hbox{  } (T^*,S^*)=\hbox{ }\hbox{ \rm ind}\hbox{  } (S,T).
$$\rm
\indent {\it Proof. \rm } Since $(T^*, S^*)^*=(S,T)$, it is enough to prove the first part of the proposition.\par
\indent First of all, note that since $S^*T^*=(TS)^*$ and $T^*S^*=(ST)^*$,
then, if $(S,T)\in P(H_1,H_2)$, $R(S^*T^*)$ and $R(T^*S^*)$
are finite dimensional subspaces of $H_1$ and $H_2$ respectively.\par
\indent Now well, according to Theorem 4.2, if $(S,T)\in P(H_1,H_2)$,
then $S^*+T\in L(H_2,H_1)$ is a Fredholm operator. Consequently,
according to Theorem 4.2 again, $(S^*, T^*)\in P(H_2, H_1)$ and
$$
\hbox{ ind }(S,T)\hbox{  } = \hbox{  }- \hbox{ ind  }(S^*,T^*).
$$
\indent However, according to Remark 2.2 or to
the observation that follows [1, Definition 1.1], $(T^*,S^*)$ is a Fredholm
pair and      
$$
\hbox{ ind }(S,T)\hbox{  } = \hbox{  } \hbox{ ind  }(T^*,S^*).\qed
$$

\newpage
\centerline{\bf 5. STABILITY PROPERTIES}
\vskip.5cm
\indent In this section, thanks to Theorem 4.2,
the stability properties of Fredholm pairs in Hilbert spaces are proved
in a direct way. In addition, note that the hypothesis in [1, Theorem 3.1], [1, Theorem 3.2] and [2, Theorem 4] can be weakened.\par  
\indent{THEOREM 5.1.} \it Let $H_1$ and $H_2$ be
two Hilbert spaces and consider $(S,T)\in P(H_1,H_2)$. Let
$S_1\in L(H_1,H_2)$ and $T_1\in L(H_2,H_1)$ be two operators such that $R(T_1S_1)$ and $R(S_1T_1)$ are finite dimensional subspaces of $H_1$
and $H_2$ respectively. Then, there is an $\epsilon>0$ such that if
$\parallel S-S_1\parallel<\epsilon$ and $\parallel T-T_1\parallel <\epsilon$, $(S_1,T_1)$ is a Fredholm pair. Futhermore,
$$
\hbox{  \rm ind}\hbox{ } (S,T)=\hbox{ \rm ind}\hbox{  }(S_1,T_1).
$$\rm
\indent {\it Proof. \rm } First of all, according to Theorem 4.2, $S+T^*\in L(H_1,H_2)$ 
is a Fredholm operator and ind $(S,T)=$ ind $(S+T^*)$.\par 
\indent Now well, since
$$
\parallel (S+T^*) -(S_1+T_1^*)\parallel <\parallel S-S_1\parallel +\parallel T-T_1\parallel,
$$
there is an $\epsilon>0$ such that if $\parallel S-S_1\parallel<\epsilon $ and $\parallel T-T_1\parallel<\epsilon$, then $S_1+T^*_1$ is a Fredholm operator and
ind $(S+T^*)$ = ind $(S_1+T_1^*)$.\par
\indent However, according to Theorem 4.2, $(S_1,T_1)$ is a Fredholm 
pair and 
$$
\hbox{ \rm ind } (S,T)=\hbox{ \rm ind }(S_1,T_1).\qed
$$
\indent {THEOREM 5.2.} \it Let $H_1$ and $H_2$ be
two Hilbert spaces and consider $(S,T)\in P(H_1,H_2)$. Let
$K\in K(H_1,H_2)$ and $K'\in K(H_2,H_1)$ be two compact operators and 
consider $S_1=S+K$ and $T_1=T+K'$. Suppose that $R(T_1S_1)$ and $R(S_1T_1)$ are finite dimensional subspaces of $H_1$
and $H_2$ respectively. Then, $(S_1,T_1)$ is a Fredholm pair. Futhermore,
$$
\hbox{ \rm ind}\hbox{ } (S,T)=\hbox{  \rm ind}\hbox{  }(S_1,T_1).
$$\rm
\indent {\it Proof. \rm } First of all, according to Theorem 4.2, $S+T^*\in L(H_1,H_2)$ 
is a Fredholm operator and ind $(S,T)=$ ind $(S+T^*)$.\par 
\indent Now well, since
$$
S_1+T^*_1=(S+ K) +(T+K')^* = (S + T^*) + (K+K^{'*}),
$$
$S_1 +T_1^*$ is a Fredholm operator and 
ind $(S+T^*)$= ind $(S_1 +T_1^*)$. \par
\indent However, according to Theorem 4.2, $(S_1,T_1)$ is a
Fredholm pair and
$$
\hbox{ \rm ind } (S,T)=\hbox{ \rm ind }(S_1,T_1).\qed
$$

\newpage
\centerline{\bf 6. FREDHOLM CHAINS}
\vskip.5cm

\indent In this section, as an applicaton of the main results of this article, characterizations of Fredholm chains in Hilbert spaces are given. Furthermore, the stability properties and the relationship between 
adjoint operators and the objects under consideration are studied. For a
general exposition, see [6].\par

\indent In what follows, sequences of spaces and maps $(H,\delta)$ and 
the operators $S$ and $T$ associated to such a sequence
will be considered, see Remark 2.4. In addition, in order to lighten the text,  
the spaces $H_p$ and the maps $\delta_p$ will be indexed for p $\in\Bbb Z$.
However, recall that there is always an n $\in\Bbb N$ such that
$H_p=0$ and $\delta_p=0$, for p $<$ 0 and p $\ge$ n+1.  
\par
\indent In first place, two characterizations
of Fredholm chains are given.\par

\indent {THEOREM 6.1.} \it Let $(H,\delta)$ be a sequence of space and
maps such that dim $R(\delta_p\delta_{p+1})$ is finite dimensional for each p $\in\Bbb Z$.
Then, the following statements are equivalent:\par
\hskip.2cm i) $(H,\delta)$ is a Fredholm chain,\par
\hskip.2cm ii) $\delta_{p+1}\delta_{p+1}^*+\delta_p^*\delta_p\in L(H_p)$
is a Fredholm self-adjoint operator for each p $\in\Bbb Z$.\par\rm
\indent {\it Proof. \rm } As in Remark 2.4, consider the spaces 
$$
H_1=\bigoplus_{p=2k}H_p,\hskip1cm H_2=\bigoplus_{p=2k+1}H_p,
$$
and the maps $S\in L(H_1,H_2)$ and $T\in L(H_2,H_1)$ defined by 
$$
S=\bigoplus_{p=2k}\hbox{ } \delta_p,\hskip1cm T=\bigoplus_{p=2k+1}\hbox{ } \delta_p.
$$
\indent According to Remark 2.4 and to Theorem 3.5, $(H,\delta)$ is a Fredholm
chain if and only if $TT^*+S^*S\in L(H_1)$ and $SS^* +T^*T\in L(H_2)$ 
are Fredholm self-adjoint operators, which is equivalent to the 
second statement of the theorem.\qed\par

\indent {THEOREM 6.2.} \it Let $(H,\delta)$ be a sequence of space and
maps such that dim $R(\delta_p\delta_{p+1})$ is finite dimensional for each p $\in\Bbb Z$.
Then, the following statements are equivalent:\par
\hskip.2cm i) $(H,\delta)$ is a Fredholm chain,\par
\hskip.2cm ii) $\oplus_{p=2k}(\delta_{p}+\delta_{p+1}^*)$ is a Fredholm operator,\par
\hskip.2cm ii) $\oplus_{p=2k+1}(\delta_{p}+\delta_{p+1}^*)$ is a Fredholm operator.\par
\indent Furthermore, in this case
$$
 \hbox{\rm ind}\hbox{ } (H,\delta)\hbox{ }=\hbox{ }\hbox{\rm ind }\hbox{ }\oplus_{p=2k}(\delta_{p}+\delta_{p+1}^*)\hbox{ }
= -\hbox{ }\hbox{\rm ind}\hbox{ }\oplus_{p=2k+1}(\delta_{p}+\delta_{p+1}^*).
$$
\rm
\indent {\it Proof. \rm } As in Theorem 6.1, consider the spaces and maps $H_1$, $H_2$, $S$ and $T$ defined in Remark 2.4.\par
\indent According to Remark 2.4 and to Theorem 4.2, $(H,\delta)$
is a Fredholm chain if and only if $S+T^*$ is a Fredholm operator, 
which is equivalent to the second statement of the theorem.\par
\indent Similarly, $(H,\delta)$
is a Fredholm pair if and only if $T+S^*$ is a Fredholm pair, which is equivalent
to the third statement of the theorem.\par
\indent Finally, the index formula is a consequence of Remark 2.4 and Theorem 4.2.\qed\par
\indent Next, the dual of a Fredholm chain is considered.\par
\indent {\it Definition\rm} 6.3 Let $(H,\delta)$ be sequence of spaces and maps
$$
0\to H_n\xrightarrow{\delta_n} H_{n-1}\to\ldots\to H_1\xrightarrow{\delta_1} H_0\to 0,
$$
where $H_p$ are Hilbert spaces, and $\delta_p\in L(H_p,H_{p-1})$ are bounded operators. In addition, assume that $H_p=0$ and $\delta_p=0$, for
p $<$ 0 and p $\ge$ n+1, where n is the first natural number with this property. The dual sequence of $(H,\delta)$ is the sequence
$(H{'},\delta{'})$, where $H{'}_p=H_{n-p}$ and $\delta{'}_p=\delta^*_{n-p+1}$,
that is the sequence
$$
0\to H_0\xrightarrow{\delta^*_1} H_{1}\to\ldots\to H_{n-1}\xrightarrow{\delta^*_n} H_n\to 0.
$$
\indent {THEOREM 6.4.} \it Let $(H,\delta)$ be a sequence of spaces and maps,
and consider its dual sequence $(H{'},\delta{'})$. Then, $(H,\delta)$ is a 
Fredholm chain if and only if $(H{'},\delta{'})$ is. Futhermore, in this case
$$
\hbox{\rm ind}\hbox{ } (H,\delta) =(-1)^n \hbox{ \rm ind}\hbox{ } (H{'},\delta{'}), 
$$
where n is the first natural number such that $H_p=0$ for p $\ge$ n+1.\par 
\rm 
\indent {\it Proof. \rm } As in Remark 2.4, consider the spaces and maps defined by $(H,\delta)$, that is
$H_1$, $H_2$, $S$ and $T$. Similarly, consider the spaces and maps defined by $(H{'},\delta{'})$,
that is $H{'}_1$, $H{'}_2$, $S{'}$ and $T{'}$. \par
\indent Now well, if n is an even natural number, then, according to Remark 2.4,
$H{'}_1=H_1$, $H{'}_2=H_2$, $S{'}=T^*$ and $T{'}=S^*$. Therefore, according to Remark 2.4 and to Theorem 4.3,
$(H,\delta)$ is a Fredholm chain if and only if $(H{'},\delta{'})$ is, and in this case $\hbox{\rm ind } (H,\delta) = \hbox{\rm ind } (H{'},\delta{'})$. \par
\indent On the other hand, if n is an odd natural number, then, according to Remark 2.4,
$H{'}_1=H_2$, $H{'}_2=H_1$, $S{'}=S^*$ and $T{'}=T^*$. Consequently, according to Remarks  2.2 and 2.4,
and to Theorem 4.3, $(H,\delta)$ is a Fredholm chain if and only if $(H{'},\delta{'})$ is, and in this case $\hbox{\rm ind } (H,\delta)$ = - $\hbox{\rm ind } (H{'},\delta{'})$. \qed\par  

\indent In the following theorems, the stability properties of Fredholm
chains in Hilbert spaces are considered.\par
\indent {THEOREM 6.5.} \it Let $(H,\delta)$ be a Fredholm chain and $(H,\delta{'})$ a sequence of space and maps such that $R(\delta{'}_p\delta{'}_{p+1})$ is finite dimensional for each p $\in\Bbb Z$. Then, there is an $\epsilon>0$ such that if $\parallel
\delta_p-\delta{'}_p\parallel<\epsilon$, p $\in\Bbb Z$, then, $(H,\delta{'})$
is a Fredholm pair. Futhermore,
$$
\hbox{ \rm ind}\hbox{ } (H,\delta) = \hbox{ \rm ind}\hbox{ } (H,\delta{'}). 
$$\rm
\indent {\it Proof. \rm } As in Theorems 6.1 and 6.2, consider the spaces and maps
$H_1$, $H_2$, $S$ and $T$ associated to $(H,\delta)$ and defined in Remark 2.4. Similarly, consider the
maps $S{'}$ and $T{'}$ associated to $(H,\delta{'})$. \par 
\indent Since 
$$
\parallel S-S{'}\parallel\le\Sigma_{p=2k}\parallel \delta_p-\delta{'}_p\parallel,
\hskip.2cm 
\parallel T-T{'}\parallel\le\Sigma_{p=2k+1}\parallel \delta_p-\delta{'}_p\parallel, 
$$
according to Theorem 5.1, there is an $\epsilon>0 $ such that if $\parallel \delta_p-\delta{'}_p\parallel<\epsilon$, p $\in \Bbb Z$, then $(S{'},T{'})$ is a Fredholm pair. Therefore, according to Remark 2.4, under this assumption
$(H,\delta{'})$ is a Fredholm chain. \par
\indent Finally, the index formula is a consequence of Remark 2.4 and Theorem 5.1.\qed
\par

\indent {THEOREM 6.6.} \it Let $(H,\delta)$ be a Fredholm chain, and let $k_p$ belong to $K(H_p,H_{p-1})$, p $\in\Bbb Z$. Consider the sequences of spaces and maps $(H,\delta{'})$, where $\delta{'}_p=\delta_p+k_p$, p $\in\Bbb Z$, and
suppose that $R(\delta{'}_p\delta{'}_{p+1})$ is finite dimensional for each p $\in\Bbb Z$. Then, 
$(H,\delta{'})$ is a Fredholm chain. Furthermore, 
$$
\hbox{\rm ind}\hbox{ } (H,\delta) = \hbox{\rm ind}\hbox{ } (H,\delta{'})
$$\rm
\indent {\it Proof. \rm } As in Theorem 6.5, consider the spaces and maps
$H_1$, $H_2$, $S$, $T$, $S{'}$ and $T{'}$.\par
\indent It is clear that $S{'}=S+K$ and $T{'}=T+K{'}$, where
$$
K=\oplus_{p=2k}k_p,\hskip.5cm K{'}=\oplus_{p=2k+1}k_p.
$$
\indent Therefore, according to Theorem 5.2 and to Remark 2.4, $(H,\delta{'})$
is a Fredholm chain.\par
\indent Finally, the index formula is a consequence of Remark 2.4 and Theorem 5.2.\qed
\par
\vskip.5cm
\indent {\bf Acknowledgments.} The author wishes to express his indebtedness to Professors C.-G. Ambrozie and
F.-H. Vasilescu. These researchers sent kindly to the author several works authored by them, which have been useful for the elaboration of the present article.\par

\vskip.5cm

\noindent Enrico Boasso\par
\noindent E-mail address: enrico\_odisseo\@yahoo.it

\begin{thebibliography}{20}

\bibitem{A} C.-G. Ambrozie, On Fredholm index in Banach spaces,
 Integral Equations Operator Theory 25 (1996), 1-34.

\bibitem{A2} C.-G. Ambrozie, The Euler characteristic is stable under
compact perturbations, Proc. Amer. Math. Soc. 124  (19969, 2041-2050.

\bibitem{CV} Z. Ceau\c sescu and F.-H. Vasilescu, Tensor products and
Taylor's joint Spectrum, Studia Math. (62) (1978), 305-311.

\bibitem{CV2} Z. Ceau\c sescu and F.-H. Vasilescu, Tensor producs and  
the joint spectrum in Hilbert spaces, Proc. Amer. Math. Soc. 72  (1978), 505-508.
 	
\bibitem{GV} C. Grosu and F.-H. Vasilescu, The K\" unneth formula for
Hilbert complexes, Integral Equations Operator Theory
 5 (1982), 1-17.

\bibitem{M} V. M\" uller, Stability of index for semi-Fredholm chains,
J. Operator Theory 37 (1997), 247-261.
 
\bibitem{V} M. Putinar, Some invariants for semi-Fredholm systems
of essentially commuting operators, J. Operator Theory 8 (1982), 65-90.

\bibitem{V} F.-H. Vasilescu, A characterization of the joint spectrum in Hilbert
spaces, Rev. Roum. Math. Pures Appl. 22 (1977), 1003-1009.

\bibitem{V2} F.-H. Vasilescu, On pairs of commuting operators, Studia Math. 62 (1978), 203-207.

\bibitem{A3} F.-H. Vasilescu, Analytic perturbations of the $\delta$-operator and integral representation fromulas in Hilbert spaces, J. Operator Theory, 1 (1979), 187-205.

\bibitem{A4} F.-H. Vasilescu, The stability of the Euler characteristic for Hilbert
complexes, Math. Ann. 248 (1980), 109-116.

\bibitem{V5} F.-H. Vasilescu, Analytic functional calculus and spectral
decompositions,  Ed. Academiei-D. Riedel Co., Bucharest-Dordrecht
(1982).
\end{thebibliography}
\end{document}